\documentclass[10pt]{article}
\usepackage{mathrsfs}
\usepackage{indentfirst}
\usepackage{amssymb,amsfonts,amsmath, psfrag,eepic,colordvi,epsfig}
\topskip 
\numberwithin{equation}{section}
\newtheorem{theorem}{Theorem}[section]

\newtheorem{lemma}[theorem]{Lemma}

\begin{document}
\parskip 7pt

\pagenumbering{arabic}
\def\sof{\hfill\rule{2mm}{2mm}}
\def\ls{\leq}
\def\gs{\geq}
\def\SS{\mathcal S}
\def\qq{{\bold q}}
\def\MM{\mathcal M}
\def\TT{\mathcal T}
\def\EE{\mathcal E}
\def\lsp{\mbox{lsp}}
\def\rsp{\mbox{rsp}}
\def\pf{\noindent {\it Proof.} }
\def\mp{\mbox{pyramid}}
\def\mb{\mbox{block}}
\def\mc{\mbox{cross}}
\def\qed{\hfill \rule{4pt}{7pt}}
\def\pf{\noindent {\it Proof.} }
\textheight=22cm

\begin{center}
{\Large\bf Some identities
  on Lin-Peng-Toh's partition  statistic of
$k$-colored partitions}
\end{center}

\begin{center}

Yang Lin$^{1}$, Ernest X.W. Xia$^{2}
 $ and   Xuan Yu$^{3}$

$^{1,2,3}$School of Mathematical Sciences, \\
  Suzhou University of Science and
Technology, \\
 Suzhou,  215009, Jiangsu Province,
 P. R. China

email: $^{1}$lin95245@163.com,
 $^{2}$ernestxwxia@163.com,  $^{3}$yuxuanusts@163.com

\end{center}

\noindent {\bf Abstract.} Recently, Andrews
 proved two conjectures
  on a partition  statistic  introduced by
     Beck.
 Very recently, Chern established
    some results on
     weighted rank and crank moments and proved
      many Andrews-Beck type congruences.
       Motivated by Andrews and Chern's work,
        Lin, Peng and Toh
  introduced
 a partition  statistic
  of $k$-colored partitions
   $NB_k(r,m,n)$
    which counts
      the total number of parts of $\pi^{(1)}$
       in each $k$-colored partition $\pi$
        of
      $n$ with ${\rm crank}_k(\pi)$
       congruent to $r$ modulo $m$
    and proved a number of congruences
     for $NB_k(r,m,n)$.
      In this paper, we prove some identities
       on $NB_k(r,m,n)$
       which are analogous to
            Ramanujan's ``most beautiful
     identity". Moreover,
      those identities imply
        some congruences proved by Lin, Peng and Toh.



\section{Introduction}

\allowdisplaybreaks

Recall that a partition $\lambda$ of a positive integer $n$ is a
sequence of positive integers $\lambda_1\geq \lambda_2\geq \cdots
\geq \lambda_k>0$ such that $\lambda_1
  +\lambda_2+\cdots
  +\lambda_k=n$. The $\lambda_i$ are called the parts of the partition
   \cite{Andrews-1976}.
   As usual, let $\#(\lambda)$, $\omega(\lambda)$, and $l(\lambda)$
   denote the number of parts of $\lambda$, the number of ones in $\lambda$
   and the largest part in $\lambda$, respectively.

   In the theory of integer partition, two
 important
  partition statistics, rank and crank,
   were defined  by
 Dyson \cite{Dyson}, and Andrews and Garvan
 \cite{Andrews-Garvan}, respectively.
 In 1944,  Dyson \cite{Dyson} defined
  the rank of a partition of $n$,
which is the largest part of the partition minus the number of
parts, i.e.,
\[
{\rm rank}(\lambda):=l(\lambda)-\#(\lambda).
\]
In 1988,
   Andrews and Garvan
  \cite{Andrews-Garvan} defined
   the crank of $\lambda$   by
\begin{align*}
\operatorname{crank}(\lambda): =\left\{\begin{array}{ll}
\l(\lambda),  & \text { if }
 \omega(\lambda)
 =0,  \\
\mu(\lambda)-\omega(\lambda),
  & \text { otherwise,}
\end{array}\right.
\end{align*}
where $\mu(\lambda)$ is the number of parts larger
 than $\omega(\lambda)$.
The two partition statistics
 are attracting broad research
interest   since
 rank
   can provide combinatorial
  interpretations of
Ramanujan's famous  congruences
\begin{align}
p(5 n+4) \equiv 0 & \pmod 5,\label{v-1}\\[6pt]
p(7 n+5) \equiv 0 & \pmod 7, \label{v-2}
\end{align}
while crank
   can provide combinatorial
  interpretations of \eqref{v-1},
  \eqref{v-2} as well as
   Ramanujan's third  congruences modulo 11:
\begin{align*}
p(11 n+6) \equiv 0 & \pmod {11},
\end{align*}
where  $p(n)$
     counts  the number
      of partitions  of $n$.

Recently, Andrews \cite{Andrews} mentioned
  that George Beck
    defined two partition   statistics
     $NT (r,m,n)$
        and $M_{\omega}(r,m,n)$,
         which
    count  the total
      number of parts in partitions of $n$
       with rank congruent to
       $r$    modulo $m$
       and the total number of ones in
        the partition
  of $n$ with crank congruent to
   $r$ modulo $m$, respectively.
  Andrews \cite{Andrews} also  proved
        the following Andrews-Beck
         type congruences
         which were conjectured by Beck:
  \begin{align*}
NT(1,5,5n+i)+2NT(2,5,5n+i)
 -2NT(3,5,5n+i) -NT(4,5,5n+i)\equiv 0
 \pmod 5
  \end{align*}
for $i\in \{1,4\}$
 and
\begin{align*}
&NT(1,7,7n+j)+NT(2,7,7n+j)
 - NT(3,7,7n+j) \nonumber\\[6pt]
 &\qquad  +NT(4,7,7n+j)
  - NT(5,7,7n+j) -NT(6,7,7n+j)\equiv 0
 \pmod 7
  \end{align*}
  for $j\in \{1,5\}$.
Very recently, Chern \cite{Chern-1,Chern-2,Chern-3} proved some
 identities involving the weighted
  rank and crank moments and established
   a number of   Andrews-Beck type congruences.
   For example, Chern \cite{Chern-1} proved that for
    $n\geq 0$,
\[
\sum_{m=1}^4 mM_{\omega}(m,5,5n+4)\equiv
 0 \pmod 5.
\]

Motivated by Andrews and Chern's work,
 Lin, Peng and Toh
 \cite{Lin} considered the generalized crank defined
  by Fu and Tang \cite{Fu}
   for $k$-colored partitions,
    where $k\geq 2$.
 To describe
this, we recall that a $k$-colored
 partition   $\pi$ of
      a positive integer $n$ is
      a $k$-tuple of partitions
      $\pi:=(\pi^{(1)},\pi^{(2)},\ldots,
       \pi^{(k)})$
      such that $|\pi^{(1)}|+|\pi^{(2)}|
      +\cdots+|\pi^{(k)}|=n$. If $\pi$
       is a $k$-colored partition of $n$, we
        denote it by $\pi \vdash
       n$.
       Fu and Tang \cite{Fu}  defined a generalized crank
        for $k$-colored partitions by
        \[
{\rm crank}_k(\pi):=\#(\pi^{(1)})-\#(\pi^{(2)}),
        \]
where $\#(\pi^{(i)})$ denotes the number of parts in
 $\pi^{(i)}$.
 Let $m,j,n,k$ be integers with $j\geq 2$, $k\geq 2$,
 $n\geq 1$ and $0\leq m \leq j-1$. Define
  \begin{align} \label{1-1}
NB_k(m,j,n):=\sum_{\pi\vdash n,
 \atop
 {\rm   crank}_k(\pi)\equiv m \
 ({\rm mod\  j})}
  \#(\pi^{(1)}).
  \end{align}
 Lin, Peng
     and Toh \cite{Lin} proved many Andrews-Beck
      type congruences for  $NB_k(m,j,n)$. For example,
       they proved that
       for $n\geq 0$,
\begin{align} \label{1-2}
\sum_{m=1}^4 mNB_2(m,5,5n+i) &
  \equiv 0 \pmod 5,\\ \label{1-3}
\sum_{m=1}^4 mNB_3(m,5,5n ) &
  \equiv 0 \pmod 5,\\ \label{1-4}
  \sum_{m=1}^4 mNB_4(m,5,5n+j) &
  \equiv 0 \pmod 5,\\ \label{1-5}
  \sum_{m=1}^4 m^3
   NB_5(m,5,5n+t) &
  \equiv 0 \pmod 5,
\end{align}
where $i\in\{0,2,3,4\}$,
 $j\in \{0,3,4\}$
  and $t\in\{2,4\}$. At the end of
   their paper, Lin, Peng
     and Toh also posed several conjectures
      on congruences for $NB_k(r,m,n)$
       which were confirmed by Du and Tang
       \cite{Du,Du-1},
        and Yao  \cite{Yao}. For more details on congruences
         for  Beck's  partition
 statistics, see
 \cite{Chan,Chen,mao,mao-1,Mao-Xia}.

In this paper, we establish some identities
 on $NB_k(r,m,n)$
  which are analogous to
 Ramanujan's ``most beautiful identity"
   \[
\sum_{n=0}^\infty p(5n+4)q^n=5\frac{(q^5;q^5)_\infty^5
 }{(q;q)_\infty^6},\]
where here and throughout this paper, we adopt the standard
  $q$-series notation
  \[
(a;q)_\infty=\prod_{n=0}^\infty (1-aq^n)
  \]
and for each positive integer $t$, set
\[
(a_1,a_2,\ldots,a_t;q)_\infty
 =(a_1;q)_\infty(a_2;q)_\infty
  \cdots (a_t;q)_\infty.
\]
Furthermore, those identities imply
 congruences \eqref{1-2}--\eqref{1-5}.

The  main results of this paper can be stated
 as follows.

\begin{theorem} \label{Th1-1}
We have
\begin{flalign}
&\sum_{n\geq0}[NB_{2}(1,5,5n)-3NB_{2}
(2,5,5n)+3NB_{2}(3,5,5n)-NB_{2}(4,5,5n)]q^n
=5q\frac{(q,q^4,q^5,q^5;q^5)_{\infty}} {(q^2,q^3;q^5)_{\infty}^3},
 \label{1-6}
\\[6pt]
&\sum_{n\geq0}[NB_{2}(1,5,5n+2) +2NB_{2}(2,5,5n+2)-2NB_{2}(3,5,5n+2)
-NB_{2}(4,5,5n+2)]q^n\nonumber\\[6pt]
&\qquad \qquad
=5\frac{(q^5;q^5)_{\infty}^2}{(q,q^4;q^5)_{\infty}^2},
 \label{1-7}
\\[6pt]
& NB_{2}(1,5,5n+3)+2NB_{2}(2,5,5n+3)
-2NB_{2}(3,5,5n+3)-NB_{2}(4,5,5n+3)  =0,
 \label{1-8}
\\[6pt]
  &\sum_{n\geq0}[NB_{2}(1,5,25n+23)-NB_{2}(4,5,25n+23)]q^n
=10\frac{(q^5;q^5)_{\infty}^5}{(q;q)_{\infty}^3},
 \label{1-9}
\\[6pt]
&\sum_{n\geq0}[NB_{2}(2,5,25n+23)-NB_{2}(3,5,25n+23)]q^n
=-5\frac{(q^5;q^5)_{\infty}^5 }{(q;q)_{\infty}^3},
 \label{1-10}
\\[6pt]
&\sum_{n\geq0}[NB_{2}(1,5,5n+4) +2NB_{2}(2,5,5n+4)
-2NB_{2}(3,5,5n+4)-NB_{2}(4,5,5n+4)]q^n
\nonumber\\[6pt]
&\qquad \qquad =5\frac{(q^5;q^5)_{\infty}^2}{(q^2,q^3
;q^5)_{\infty}^2}. \label{1-11}
\end{flalign}
\end{theorem}

\begin{theorem} \label{Th1-2}
We have
\begin{flalign}
&\sum_{n\geq0}[NB_{3}(1,5,5n)-NB_{3}(4,5,5n)+2NB_{3}(2,5,5n)-2NB_{3}(3,5,5n)]q^n
\nonumber\\[6pt]
=&\frac{(q^2,q^3; q^5)_{\infty}(q^5;q^5)_{\infty}
}{2(q,q^4;q^5)_{\infty}^3} -\frac{(q^2,q^3;q^5)_{\infty}(
q^5;q^5)_{\infty}^2}{2(q;q)_{\infty}^5
(q,q^4;q^5)_{\infty}^3}-5\frac{(q^2,q^3;
q^5)_{\infty}(q^5;q^5)_{\infty}^2
}{(q;q)_{\infty}^6(q,q^4;q^5)_{\infty}^2}
\sum_{n=-\infty}^{\infty}(-1)^n n q^{\frac{n(5n+1)}{2}}
\nonumber\\[6pt]
&+\frac{15q(q^2,q^3;q^5)_{\infty}
(q^5;q^5)_{\infty}^7}{2(q;q)_{\infty}^6
(q,q^4;q^5)_{\infty}^3}\bigg(4 +3q\frac{(q,q^4;q^5)_{\infty}^5
}{(q^2,q^3;q^5)_{\infty}^5}\bigg) -\frac{5q(q^5;q^5)_{\infty}^5}{2
(q;q)_{\infty}^6(q^2,q^3;q^5)_{\infty}}. \label{1-12}
\end{flalign}
\end{theorem}

\begin{theorem} \label{Th1-3}
We have
\begin{flalign}
&\sum_{n\geq 0}[NB_4(1,5,5n)-NB_4(4,5,5n)+2NB_4(2,5,5n)-2NB_4(3,5,5
n)] q^n
  \nonumber \\[6pt]
=&5\sum_{n=-\infty}^{\infty}(-1)^nnq^{\frac{n(5n+1)}{2}}
\bigg(-7q^3\frac{(q,q^4;q^5)_{\infty}^6(q^5;q^5)_{\infty}^9}{(q^2,q^3;q^5)_{\infty}^8(q;q)_{\infty}^{12}}
+12q^2\frac{(q,q^4;q^5)_{\infty}(q^5;q^5)_{\infty}^9}{(q^2,q^3;q^5)_{\infty}^3(q;q)_{\infty}^{12}}
  \nonumber \\[6pt]
&-4q\frac{(q^2,q^3;q^5)_{\infty}^2(q^5;q^5)_{\infty}^9}{(q,q^4;q^5)_{\infty}^4(q;q)_{\infty}^{12}}
-9\frac{(q^2,q^3;q^5)_{\infty}^7(q^5;q^5)_{\infty}^9}{(q,q^4;q^5)_{\infty}^9(q;q)_{\infty}^{12}}
+20q^2\frac{(q,q^4;q^5)_{\infty}^4(q^5;q^5)_{\infty}^7}{(q^2,q^3;q^5)_{\infty}^5(q;q)_{\infty}^{12}}
  \nonumber \\[6pt]
&-5q\frac{(q^5;q^5)_{\infty}^7}{(q,q^4;q^5)_{\infty}(q;q)_{\infty}^{12}}
+5\frac{(q^2,q^3;q^5)_{\infty}^5(q^5;q^5)_{\infty}^7}{(q,q^4;q^5)_{\infty}^6(q;q)_{\infty}^{12}}\bigg)
  \nonumber \\[6pt]
&+10q^2\frac{(q,q^4;q^5)_{\infty}^4(q^5;q^5)_{\infty}^8}{(q^2,q^3;q^5)_{\infty}^4(q;q)_{\infty}^{12}}
-\frac{5}{2}q\frac{(q^2,q^3;q^5)_{\infty}(q^5;q^5)_{\infty}^8}{(q,q^4;q^5)_{\infty}(q;q)_{\infty}^{12}}
+\frac{5}{2}\frac{(q^2,q^3;q^5)_{\infty}^6(q^5;q^5)_{\infty}^8}{(q,q^4;q^5)_{\infty}^6(q;q)_{\infty}^{12}}
  \nonumber \\[6pt]
&-\frac{305}{2}q^3\frac{(q,q^4;q^5)_{\infty}^3(q^5;q^5)_{\infty}^{12}}{(q^2,q^3;q^5)_{\infty}^5(q;q)_{\infty}^{12}}
+\frac{335}{2}q^2\frac{(q^5;q^5)_{\infty}^{12}}{(q,q^4;q^5)_{\infty}^2(q;q)_{\infty}^{12}}
+85q\frac{(q^2,q^3;q^5)_{\infty}^5(q^5;q^5)_{\infty}^{12}}{(q,q^4;q^5)_{\infty}^7(q;q)_{\infty}^{12}}
  \nonumber \\[6pt]
&+5\frac{(q^2,q^3;q^5)_{\infty}^{10}(q^5;q^5)_{\infty}^{12}}{(q,q^4;q^5)_{\infty}^{12}(q;q)_{\infty}^{12}}
+25q^3\frac{(q,q^4;q^5)_{\infty}^6(q^5;q^5)_{\infty}^{10}}{(q^2,q^3;q^5)_{\infty}^7(q;q)_{\infty}^{12}}
-15q^2\frac{(q,q^4;q^5)_{\infty}(q^5;q^5)_{\infty}^{10}}{(q^2,q^3;q^5)_{\infty}^2(q;q)_{\infty}^{12}}
  \nonumber \\[6pt]
&-10q\frac{(q^2,q^3;q^5)_{\infty}^3(q^5;q^5)_{\infty}^{10}}{(q, q^4
; q^5)_{\infty}^4(q;q)_{\infty}^{12}}
-\frac{15}{2}\frac{(q^2,q^3;q^5)_{\infty}^8(q^5;q^5)_{\infty}^{10}}{(q,q^4;
q^5)_{\infty}^9(q;q)_{\infty}^{12}}
 , \label{1-13} \\[6pt]
&\sum_{n\geq 0}[N B_4(1,5,5 n+3)-N B_4(4,5,5 n+3)+2 N B_4(2,5,5
n+3)-2 N B_4(3,5,5 n+3)] q^n
  \nonumber \\[6pt]
=&55q^3 \frac{(q, q^4 ; q^5)_{\infty}^6(q^5
;q^5)_{\infty}^{12}}{(q^2, q^3 ;
q^5)_{\infty}^8(q;q)_{\infty}^{12}}-\frac{335}{2} q^2
\frac{(q,q^4;q^5)_{\infty}(q^5 ; q^5)_{\infty}^{12}}{(q^2, q^3;
q^5)_{\infty}^3(q;q)_{\infty}^{12}} +260 q \frac{(q^2, q^3 ;
q^5)_{\infty}^2(q^5 ; q^5)_{\infty}^{12}}{(q, q^4 ;
q^5)_{\infty}^4(q;q)_{\infty}^{12}}
  \nonumber \\[6pt]
&+\frac{45}{2} \frac{(q^2, q^3 ; q^5)_{\infty}^{7}(q^5 ;
q^5)_{\infty}^{12}}{(q, q^4 ; q^5)_{\infty}^{9}(q;q)_{\infty}^{12}}
 -5q^2
\frac{(q, q^4 ; q^5)_{\infty}^4(q^5;q^5)_{\infty}^{10}}{(q^2, q^3 ;
q^5)_{\infty}^5(q;q)_{\infty}^{12}}
 -\frac{35}{2}q \frac{(q^5;q^5)_{\infty}^{10}}{(q, q^4 ;
q^5)_{\infty}(q;q)_{\infty}^{12}}
  \nonumber \\[6pt]
&-\frac{15}{2} \frac{(q^2, q^3 ;
q^5)_{\infty}^5(q^5;q^5)_{\infty}^{10}}{(q, q^4 ;
q^5)_{\infty}^6(q;q)_{\infty}^{12}}
, \label{1-14} \\[6pt]
&\sum_{n\geq 0}[N B_4(1,5,5 n+4)-N B_4(4,5,5 n+4)+2 N B_4(2,5,5
n+4)-2 N B_4(3,5,5 n+4)] q^n
  \nonumber \\[6pt]
=&-\frac{45}{2}q^3\frac{(q,q^4;q^5)_{\infty}^7(q^5;q^5)_{\infty}^{12}}{(q^2,q^3;q^5)_{\infty}^9(q;q)_{\infty}^{12}}
+260q^2\frac{(q,q^4;q^5)_{\infty}^2(q^5;q^5)_{\infty}^{12}}{(q^2,q^3;q^5)_{\infty}^4(q;q)_{\infty}^{12}}
+\frac{335}{2}q\frac{(q^2,q^3;q^5)_{\infty}(q^5;q^5)_{\infty}^{12}}{(q,q^4;q^5)_{\infty}^3(q;q)_{\infty}^{12}}
  \nonumber \\[6pt]
&+55\frac{(q^2,q^3;q^5)_{\infty}^{6}(q^5;q^5)_{\infty}^{12}}{(q,q^4;q^5)_{\infty}^{8}(q;q)_{\infty}^{12}}
+\frac{5}{2}q^2\frac{(q,q^4;q^5)_{\infty}^5(q^5;q^5)_{\infty}^{10}}{(q^2,q^3;q^5)_{\infty}^6(q;q)_{\infty}^{12}}
+\frac{5}{2}q
\frac{(q^5;q^5)_{\infty}^{10}}{(q^2,q^3;q^5)_{\infty}(q;q)_{\infty}^{12}}
  \nonumber \\[6pt]
&-15\frac{(q^2,q^3;q^5)_{\infty}^4
(q^5;q^5)_{\infty}^{10}}{(q,q^4;q^5
)_{\infty}^5(q;q)_{\infty}^{12}}. \label{1-15}
\end{flalign}
\end{theorem}

\begin{theorem} \label{Th1-4}
We have
\begin{flalign}
&\sum_{n\geq 0}[N B_5(1,5,5 n+2)-N B_5(4,5,5 n+2)+3 N B_5(2,5,5
n+2)-3N B_5(3,5,5 n+2)] q^n
  \nonumber \\[6pt]
=&-70q^5\frac{(q,q^4;q^5)_{\infty}^{11}(q^5;q^5)_{\infty}^{17}}{(q^2,q^3;q^5)_{\infty}^{13}(q;q)_{\infty}^{18}}
+2165q^4\frac{(q,q^4;q^5)_{\infty}^{6}(q^5;q^5)_{\infty}^{17}}{(q^2,q^3;q^5)_{\infty}^{8}(q;q)_{\infty}^{18}}
-\frac{3265}{2}q^3 \frac{(q, q^4 ; q^5)_{\infty}(q^5
;q^5)_{\infty}^{17}}{(q^2, q^3 ; q^5)_{\infty}^3(q;q)_{\infty}^{18}}
  \nonumber \\[6pt]
&+\frac{5745}{2}q^2 \frac{(q^2,q^3;q^5)_{\infty}^2(q^5
;q^5)_{\infty}^{17}}{(q, q^4; q^5)_{\infty}^4(q;q)_{\infty}^{18}}
 +\frac{1965}{2}q
\frac{(q^2, q^3 ; q^5)_{\infty}^7(q^5 ; q^5)_{\infty}^{17}}{(q, q^4
; q^5)_{\infty}^9(q;q)_{\infty}^{18}} +\frac{35}{2}\frac{(q^2, q^3 ;
q^5)_{\infty}^{12}(q^5 ; q^5)_{\infty}^{17}}{(q, q^4 ;
q^5)_{\infty}^{14}(q;q)_{\infty}^{18}}
  \nonumber \\[6pt]
&+10q^4\frac{(q,q^4;q^5)_{\infty}^9(q^5;q^5)_{\infty}^{15}}{(q^2,q^3;q^5)_{\infty}^{10}(q;q)_{\infty}^{18}}
-45q^3\frac{(q,q^4;q^5)_{\infty}^4(q^5;q^5)_{\infty}^{15}}{(q^2,q^3;q^5)_{\infty}^5(q;q)_{\infty}^{18}}
-\frac{165}{2}q^2 \frac{(q^5;q^5)_{\infty}^{15}}{(q, q^4 ;
q^5)_{\infty}(q;q)_{\infty}^{18}}
  \nonumber \\[6pt]
&-200q \frac{(q^2,q^3;q^5)_{\infty}^5(q^5;q^5)_{\infty}^{15}}{(q,
q^4 ; q^5)_{\infty}^6(q;q)_{\infty}^{18}} -\frac{15}{2} \frac{(q^2,
q^3 ; q^5)_{\infty}^{10}(q^5;q^5)_{\infty}^{15}}{(q, q^4 ;
q^5)_{\infty}^{11}(q;q)_{\infty}^{18}}
  , \label{1-16} \\[6pt]
&\sum_{n\geq0}[NB_5(1,5,5n+4)-NB_5(4,5,5n+4)+3NB_5(2,5,5n+4)-3NB_5(3,5,5n+4)]q^n
  \nonumber \\[6pt]
=&5\sum_{n=-{\infty}}^\infty(-1)^nnq^{n(5n+1)/2}\bigg(
2q^4\frac{(q,q^4;q^5)_{\infty}^{11}(q^5;q^5)_{\infty}^{14}}{(q^2,q^3;q^5)_{\infty}^{13}(q;q)_{\infty}^{18}}
-q^3\frac{(q,q^4;q^5)_{\infty}^{6}(q^5;q^5)_{\infty}^{14}}{(q^2,q^3;q^5)_{\infty}^{8}(q;q)_{\infty}^{18}}
  \nonumber \\[6pt]
&+2108q^2\frac{(q,q^4;q^5)_{\infty}(q^5;q^5)_{\infty}^{14}}{(q^2,q^3;q^5)_{\infty}^3(q;q)_{\infty}^{18}}
-4q\frac{(q^2,q^3;q^5)_{\infty}^{2}(q^5;q^5)_{\infty}^{14}}{(q,q^4;q^5)_{\infty}^{4}(q;q)_{\infty}^{18}}
-3\frac{(q^2,q^3;q^5)_{\infty}^{7}(q^5;q^5)_{\infty}^{14}}{(q,q^4;q^5)_{\infty}^{9}(q;q)_{\infty}^{18}}
  \nonumber \\[6pt]
&-5q^3\frac{(q,q^4;q^5)_{\infty}^{9}(q^5;q^5)_{\infty}^{12}}{(q^2,q^3;q^5)_{\infty}^{10}(q;q)_{\infty}^{18}}
+30q^2\frac{(q,q^4;q^5)_{\infty}^{4}(q^5;q^5)_{\infty}^{12}}{(q^2,q^3;q^5)_{\infty}^{5}(q;q)_{\infty}^{18}}\bigg)
  \nonumber \\[6pt]
&-\frac{5}{2}q^5\frac{(q,q^4;q^5)_{\infty}^{13}(q^5;q^5)_{\infty}^{17}}{(q^2,q^3;q^5)_{\infty}^{15}(q;q)_{\infty}^{18}}
+825q^4\frac{(q,q^4;q^5)_{\infty}^{8}(q^5;q^5)_{\infty}^{17}}{(q^2,q^3;q^5)_{\infty}^{10}(q;q)_{\infty}^{18}}
+\frac{1865}{2}q^2\frac{(q^5;q^5)_{\infty
}^{17}}{(q,q^4;q^5)_{\infty}^{2}(q;q)_{\infty}^{18}}
  \nonumber \\[6pt]
&+\frac{4125}{2}q\frac{(q^2,q^3;q^5)_{\infty}^{5}(q^5;q^5)_{\infty}^{17}}{(q,q^4;q^5)_{\infty}^{7}(q;q)_{\infty}^{18}}
+\frac{275}{2}\frac{(q^2,q^3;q^5)_{\infty}^{10}(q^5;q^5)_{\infty}^{17}}{(q,q^4;q^5)_{\infty}^{12}(q;q)_{\infty}^{18}}
-5q^4\frac{(q,q^4;q^5)_{\infty}^{11}(q^5;q^5)_{\infty}^{15}}{(q^2,q^3;q^5)_{\infty}^{12}(q;q)_{\infty}^{18}}
  \nonumber \\[6pt]
&+\frac{1895}{2}q^2\frac{(q,q^4;q^5)_{\infty}(q^5;q^5)_{\infty}^{15}}{(q^2,q^3;q^5)_{\infty}^{2}(q;q)_{\infty}^{18}}
-245q\frac{(q^2,q^3;q^5)_{\infty}^3(q^5;q^5)_{\infty}^{15}}{(q,q^4;q^5)_{\infty}^4(q;q)_{\infty}^{18}}
-\frac{95}{2}\frac{(q^2,q^3;q^5)_{\infty}^8(q^5;q^5)_{\infty}^{15}}{(q,q^4;q^5)_{\infty}^9(q;q)_{\infty}^{18}}
  \nonumber \\[6pt]
&-\frac{5}{2}q^3\frac{(q,q^4;q^5)_{\infty}^9(q^5;q^5)_{\infty}^{13}}{(q^2,q^3;q^5)_{\infty}^9(q;q)_{\infty}^{18}}
+15q^2\frac{(q,q^4;q^5)_{\infty}^4
(q^5;q^5)_{\infty}^{13}}{(q^2,q^3;
q^5)_{\infty}^4(q;q)_{\infty}^{18}}. \label{1-17}
\end{flalign}
\end{theorem}

Note that Theorems \ref{Th1-1}--\ref{Th1-4}
 imply
 \eqref{1-2}--\eqref{1-5}, respectively.

 The paper is organized as follows.
 In Section 2,  we
    prove some identities  which are used
    to
     prove the main results of this paper.
      Sections  3--6  are   devoted to
         the proofs of Theorems
         \ref{Th1-1}--\ref{Th1-4}.
          In Section 7, we make some
          concluding remarks
          concerning future
 work.

\section{Preliminaries}

In order
 to prove the main results
  of this paper, we need
  some lemmas.

The following identities were
 proved by Jin, Liu and Xia \cite{Liuxin-jin}.

\begin{lemma}
 \cite[(2.14), Lemma 2.3]{Liuxin-jin}\label{L-1}
 Define
  \begin{align} \label{2-1}
X(a,b,c):= &\sum_{n=0}^{\infty} \frac{q^{bn+c}
}{1-q^{5n+a}}-\sum_{n=0}^{\infty} \frac{q^{
(5-b)n+(5+c-a-b)}}{1-q^{5n+(5-a)}},\\[6pt]
Y(d):=&\sum_{n=1}^{\infty} \frac{ q^{dn}}{1-q^{5n}}
-\sum_{n=1}^{\infty} \frac{q^{(5-d)n}
 }{1-q^{5n}},
\end{align}
where $a,b,c $ are integer with $1\leq a,b,d
 \leq 4$. Then
\begin{align}
&X(4,1,0)=0,
  \label{2-3}\\[6pt]
&X(3,1,0)=\frac{(q^5,q^5;q^5 )_{\infty}}{(q^2,q^3;q^5)_{\infty}},
  \label{2-4} \\[6pt]
&X(2,1,0)=\frac{(q^5,q^5;q^5)_{\infty} }{(q,q^4;q^5)_{\infty}},
  \label{2-5} \\[6pt]
&X(1,1,0)=\frac{(q^2,q^3,q^5,q^5;q^5)_{\infty}}{(q,q^4;q^5)_{\infty}^{2}},
  \label{2-6} \\[6pt]
&X(2,2,0)=\frac{(q,q^4,q^5,q^5;q^5)_{\infty}}{(q^2,q^3;q^5)_{\infty}^2}
 , \label{2-7} \\
&X(4,2,1)=-\frac{(q^5,q^5;q^5)_{\infty}}{(q^2,q^3;q^5)_{\infty}}
 , \label{2-8} \\[6pt]
&X(1,2,0)=\frac{(q^5,q^5;q^5)_{\infty}} {(q,q^4;q^5)_{\infty}},
\label{2-9}
\\[6pt]
&X(3,2,1)=0
 , \label{2-10}
\\[6pt]
&Y(1)=\frac{3(q^2,q^3,q^5;q^5)_{\infty}^2}{10(q,q^4;q^5)_{\infty}^3}+\frac{q(q,q^4,q^5;q^5)_{\infty}^2}{10(q^2,q^3;q^5)_{\infty}^3}
-\frac{3}{10}
 , \label{2-11} \\[6pt]
& Y(2)=\frac{(q^2,q^3,q^5;q^5 )_{\infty}^2}{10(q,q^4;q^5)_{\infty}^3
}-\frac{3q(q,q^4,q^5;q^5)_{\infty}^2
}{10(q^2,q^3;q^5)_{\infty}^3}-\frac{1}{10}. \label{2-12}
\end{align}
\end{lemma}

\begin{lemma} \label{L-2} We have
\begin{align}
&\frac{1}{R(q)^5}-11-R(q)^5=\frac{(q;q)_{\infty}^6}{q(q^5;q^5)_{\infty}^6},
  \label{2-13} \\[6pt]
&\frac{q^{\frac{3}{5}}}{R(q)^3} -3q^{\frac{3}{5}}R(q)^2 =\frac{10
}{(q^5;q^5)_{\infty}^3}\sum_{n= -\infty}^{\infty}(-1)^nn
q^{\frac{n(5n+1)}{2}}+\frac{(q^2,q^3;q^5)_\infty
}{(q^5;q^5)_{\infty}^2}
, \label{2-14}\\[6pt]
&\frac{1}{(q;q)_\infty}= \frac{(q^{25};q^{25})_\infty^5
 }{(q^5;q^5)_\infty^6}\left(
T^4 +qT^3 +2q^2T^2
 +3q^3T +5q^4-3\frac{q^5}{T }
  +2\frac{q^6}{T^2 }-\frac{q^7}{T^3 }
   +\frac{q^8}{T^4  }\right), \label{2-15}
\end{align}
where $R(q)$   is the Rogers-Ramanujan
 continued fraction
defined by
\[
R(q)=q^{\frac{1}{5}}\frac{(q,q^4;
q^5)_{\infty}}{(q^2,q^3;q^5)_{\infty}}
\]
and $ T:=\frac{V(q^5)}{W(q^5)}  $ with
\begin{align}\label{2-16}
V:=V(q)=\frac{1}{(q,q^4;q^5)_{\infty}} ,\qquad \qquad
W:=W(q)=\frac{1}{(q^2,q^3;q^5)_{\infty}}.
\end{align}
\end{lemma}

\noindent{\it Proof.}
 Identity \eqref{2-13}
   follows  from (1.1.11)
     in \cite{Andrews-0} and
      \eqref{2-15}
     follows from (7.4.14)
      in \cite{Berndt-2006}.
       The following identity was proved by Andrews
        and Berndt \cite[(1.3.2), p. 19]{Andrews-0}:
\begin{align}\label{2-17}
\frac{q^{\frac{3}{5}}}{R(q)^3}-3q^{\frac{3}{5}}R(q)^2=\frac{1}{(q^5;q^5)_{\infty}^3}\sum_{-\infty}^{\infty}(-1)^n(10n+1)q^{\frac{n(5n+1)}{2}}
.
\end{align}
By the Jacobi triple product identity,
\begin{align}\label{2-18}
\sum_{n=-\infty}^{\infty}(-1)^n q^{\frac{n(5n+1)}{2}
}=(q^2,q^3,q^5;q^5)_{\infty}.
\end{align}
Combining \eqref{2-17}
 and \eqref{2-18}, we arrive at \eqref{2-14}.
  The proof of Lemma \ref{L-2}
   is complete.
\qed

\begin{lemma}\label{L-3}
We have
\begin{align}\label{2-19}
\sum_{n\geq 0}NB_k(r,5,n)q^n =&\frac{1}{5}
 \sum_{j=0}^{4}
 \zeta_5^{-rj}\frac{(q;q)_\infty^{2-k}
  }{(\zeta_5^j q;q)_\infty
  (q/\zeta_5^j;q)_\infty}
  \left(\sum_{n=1}^\infty
   \frac{q^{5n}}{1-q^{5n}}+\sum_{t=1}^4
    \zeta_5^{tj}
    R_t(q)
   \right),
\end{align}
where $R_t(q)$ is defined by
\begin{align}\label{2-20}
R_t(q)=\sum_{n=1}^\infty
   \frac{q^{ tn}}{1-q^{5n}}.
\end{align}
\end{lemma}

\noindent{\it Proof.}
 Lin, Peng and Toh \cite[(3.2)]{Lin}
  proved that
  \begin{align}\label{2-21}
\sum_{n=0}^\infty \sum_{\pi \vdash
 n} x^{\#(\pi^{(1)})}z^{{\rm crank}_k (\pi)
  } q^n=\frac{(q;q)_\infty^{2-k}}{
   (qxz;q)_\infty (q/z;q)_\infty}.
  \end{align}
Applying  the operator $\left[
 \frac{\partial}{\partial x}\right]_{x=1}$
 to \eqref{2-21} yields
 \begin{align}
\sum_{n=0}^\infty \sum_{\pi \vdash
 n} \#(\pi^{(1)})  z^{{\rm crank}_k (\pi)
  } q^n=&\left[\frac{\partial}{\partial x}\frac{(q;q)_\infty^{2-k}}{
   (qxz;q)_\infty (q/z;q)_\infty}
   \right]_{x=1}
   \nonumber\\[6pt]
   =&\left[\frac{(q;q)_\infty^{2-k}}{
   (qxz;q)_\infty (q/z;q)_\infty}
   \frac{\partial}{\partial x} \log \left(
   \frac{1}{(qzx;q)_\infty}\right)
   \right]_{x=1}
   \nonumber\\[6pt]
   =&\frac{(q;q)_\infty^{2-k}}{
   (q z;q)_\infty (q/z;q)_\infty}
   \sum_{n=1}^\infty
    \frac{zq^n}{1-zq^n}.\label{2-22}
 \end{align}
Let $r$ and $m$ be  integers
 with $m\geq 2$ and $0\leq r \leq m-1$.
  Throughout this paper, we always let
$  \zeta_m={e^{2\pi i/m}}$.
     By \eqref{1-1},
      \eqref{2-22} and the fact
that
\begin{align*}
  \sum_{j=0}^{m-1}
\zeta_m^{kj}= &  \left\{
  \begin{aligned} & m
    ,\qquad {\rm if\ }
k
  \equiv 0 \ ({\rm  mod }\ m), \\[6pt]
 & 0
    ,\qquad \ {\rm if\ }
k
  \not\equiv 0 \   ({\rm  mod }\ m),
  \end{aligned} \right.
\end{align*}
we deduce that
 for any integer $b$,
 \begin{align}
\sum_{n\geq 0}NB_k(r,m,n)q^n =&
  \sum_{n=0}^\infty
 \sum_{\lambda\vdash n,\atop
  {\rm crank}_k(\lambda)\equiv r \pmod m}
   \#(\pi^{(1)}) q^n
\nonumber\\[6pt]
=& \sum_{n=0}^\infty
 \sum_{\lambda\vdash n} \left(
 \frac{1}{m}\sum_{j=0}^{m-1}
  \zeta_m^{({\rm crank}_k(\lambda)-r)j}\right)
   \#(\pi^{(1)}) q^n
\nonumber\\[6pt]
=&  \frac{1}{m}\sum_{j=0}^{m-1}
 \zeta_m^{-rj}
\sum_{n=0}^\infty
 \sum_{\lambda\vdash n}
  \zeta_m^{ {\rm crank}_k(\lambda) j}
  \#(\pi^{(1)}) q^n
\nonumber\\[6pt]
 = &
\frac{1}{m} \sum_{j=0}^{m-1}
 \zeta_m^{-rj}\frac{(q;q)_\infty^{2-k}
  }{(\zeta_m^j q;q)_\infty
  (q/\zeta_m^j;q)_\infty}
    \sum_{n=1}^\infty
    \frac{\zeta_m^{j}q^n}{1-q^n\zeta_m^{j}}
     . \label{2-23}
 \end{align}
Moreover, using the fact that
 \begin{align*}
(1-q^{5n})=(1-\zeta_5^{j}q^n) (1 +\zeta_5^{j}q^n +\zeta_5^{2j}q^{2n}
+\zeta_5^{3j}q^{3n} +\zeta_5^{4j}q^{4n}),
\end{align*}
 we deduce  that
\begin{align}\label{2-24}
\sum_{n=1}^\infty
 \frac{\zeta_5^{ j} q^n}{1-\zeta_5^{ j}q^{n}}
  =  \sum_{n=1}^\infty
   \frac{q^{5n}}{1-q^{5n}}+\sum_{t=1}^4
    \zeta_5^{tj}
    R_t(q)
   ,
\end{align}
where $R_t(q)$ is defined by \eqref{2-20}.
  Setting $m=5$ in \eqref{2-23}
    and using \eqref{2-24}, we
arrive at \eqref{2-19}.
 This completes the proof of Lemma \ref{L-3}.
  \qed

\section{Proof of Theorem  \ref{Th1-1}}

In \cite{Garvan}, Garvan proved that
\begin{align}\label{3-1}
\frac{1
  }{(\zeta_5 q;q)_\infty
  (q/\zeta_5;q)_\infty}
  =\frac{1
  }{(\zeta_5^4 q;q)_\infty
  (q/\zeta_5^4;q)_\infty}
  =V(q^5)+\left(\zeta_5+\frac{1}{\zeta_5 }\right)
   qW(q^5)
\end{align}
and
\begin{align}\label{3-2}
\frac{1
  }{(\zeta_5^2 q;q)_\infty
  (q/\zeta_5^2;q)_\infty}
  =\frac{1
  }{(\zeta_5^3 q;q)_\infty
  (q/\zeta_5^3;q)_\infty}
  =V(q^5)+\left(\zeta_5^2+\frac{1}{\zeta_5^2}\right)
   qW(q^5),
\end{align}
where $V(q)$ and $W(q)$ are defined by
  \eqref{2-16}.
 Setting $k=2$ in \eqref{2-19} and
applying  \eqref{3-1} and \eqref{3-2}, we deduce
 that
\begin{align} \label{3-3}
&\sum_{n=0}^{\infty}(NB_2(1,5,n)-NB_2(4,5,n))q^n=V(q^5)(R_1(q)-R_4(q))+qW(q^5)(R_2(q)-R_3(q))
\end{align}
and
\begin{align}\label{3-4}
 \sum_{n=0}^{\infty}(NB_2(2,5,n)
 -NB_2(3,5,n))q^n=&V(q^5)(R_2(q)-R_3(q))
 \nonumber\\[6pt]
 &+qW(q^5)(R_1(q)-R_2(q)+R_3(q)-R_4(q)).
\end{align}
If we extract
 those  terms in which the power of $q$
 is
  congruent to 0 modulo 5 in \eqref{3-3}, then
   replace  $q^5$ by $q$, we arrive at
\begin{flalign}
&\sum_{n=0}^{\infty}(NB_2(1,5,5n) -NB_2(4,5,5n))q^n =V(q)Y(1)+W(q)q
X(2,2,0). \label{3-5}
\end{flalign}
Substituting
 \eqref{2-7} and \eqref{2-11} into
  \eqref{3-5},
  we arrive at
\begin{flalign} \label{3-6}
&\sum_{n=0}^{\infty}(NB_2(1,5,5n)-NB_2(4,5,5n))q^n
  \nonumber \\[6pt]
 =&\frac{1}{(q,q^4;q^5)_{\infty}}\bigg(\frac{3(q^2,q^3,q^5;q^5)_{\infty}^2}{10(q,q^4;q^5)_{\infty}^3}+\frac{q(q,q^4,q^5;q^5)_{\infty}^2}{10(q^2,q^3;q^5)_{\infty}^3}-\frac{3}{10}
\bigg)+\frac{1}{(q^2,q^3;q^5)_{\infty}}\bigg(q\frac{(q,q^4,q^5,q^5;q^5)_{\infty}}{(q^2,q^3;q^5)_{\infty}^2}
\bigg)\nonumber\\[6pt]
 =&\frac{3(q^2,q^3,q^5;q^5)_{\infty}^2}
 {10(q,q^4;q^5)_{\infty}^4}+
 \frac{11q(q,q^4,q^5,q^5;q^5)_{\infty}}{
 10(q^2,q^3;q^5)_{\infty}^3}
 -\frac{3}{10(q,q^4;q^5)_{\infty}}.
\end{flalign}
Extracting
 those  terms in which the power
 of $q$ is
  congruent to 0 modulo 5 in \eqref{3-4},
   then
   replacing  $q^5$ by $q$
    and using \eqref{2-3},
      \eqref{2-7} and \eqref{2-12}, we deduce that
\begin{flalign}
&\sum_{n=0}^{\infty}(NB_2(2,5,5n)-NB_2(3,5,5n))q^n
  \nonumber\\[6pt]
 =&V(q)Y(2)
+W(q)q(X(4,1,0)-X(2,2,0))
\nonumber\\[6pt]
 =&\frac{1}{(q,q^4;q^5)_{\infty}}\bigg(\frac{(q^2,q^3,q^5;q^5)_{\infty}^2}{10(q,q^4;q^5)_{\infty}^3}-\frac{3q(q,q^4,q^5;q^5)_{\infty}^2}{10(q^2,q^3;q^5)_{\infty}^3}-\frac{1}{10}\bigg)
 -q\frac{1}{(q^2,q^3;q^5)_{\infty}}
\frac{(q,q^4,q^5,q^5;q^5)_{\infty}}{(q^2,q^3;q^5)_{\infty}^2}
\nonumber\\[6pt]
 =&\frac{(q^2,q^3,q^5;q^5
 )_{\infty}^2}{10(q,q^4;q^5)_{\infty}^4}
 -\frac{13q(q,q^4,q^5,q^5;q^5)_{\infty}
 }{10(q^2,q^3;q^5)_{\infty}^3}
 -\frac{1}{10(q,q^4;q^5)_{\infty}}.
 \label{3-7}
\end{flalign}
Identity
 \eqref{1-6} follows from \eqref{3-6}
 and \eqref{3-7}.

  Extracting those terms in which
the power of $q$ is congruent to
 $2$ modulo $5$ in \eqref{3-3} and
\eqref{3-4}, then dividing by $q^2$ and replacing $q^5$ by $q$,
  we
obtain
\begin{align}
&\sum_{n=0}^{\infty}(NB_2(1,5,5n+2)-NB_2(4,5,5n+2))q^n
  \nonumber \\[6pt]
 =&V(q) X(2,1,0) +W(q) X(3,2,1)
\nonumber \\[6pt]
 =&\frac{(q^5;q^5)_{\infty}^2}
 {(q,q^4;q^5)_{\infty}^2} \qquad
 \qquad ({\rm by}\  \eqref{2-5}\ {\rm
  and} \ \eqref{2-10}) \label{3-8}
\end{align}
and
\begin{align}
&\sum_{n=0}^{\infty}(NB_2(2,5,5n+2)-NB_2(3,5,5n+2))q^n
\nonumber\\[6pt]  =&V(q) X(1,2,0)
 +W(q) (X(1,1,0)-X(3,2,1))
\nonumber\\[6pt]  =&2\frac{(q^5;q^5) _{\infty}^2}
{(q,q^4;q^5)_{\infty}^2},  \qquad ({\rm by} \ \eqref{2-6},
\eqref{2-9}\ {\rm
  and} \ \eqref{2-10}) \label{3-9}
\end{align}
from which with
 \eqref{3-8}  and  \eqref{3-9},  identity \eqref{1-7} follows.

 Picking out those terms in which
 the power of $q$ is congruent to
$3$ modulo $5$   in \eqref{3-3} and \eqref{3-4}, then dividing by
$q^3$ and replacing $q^5$ by $q$, we can get that
\begin{align}
&\sum_{n=0}^{\infty}(NB_2(1,5,5n+3) -NB_2(4,5,5n+3))q^n=V(q)
X(3,1,0) +W(q) X(1,2,0) \label{3-10}
\end{align}
and \begin{align}
&\sum_{n=0}^{\infty}(NB_2(2,5,5n+3)-NB_2(3,5,5n+3))q^n =V(q)
X(4,2,1) +W(q) (X(2,1,0)-X(1,2,0)).
 \label{3-11}
\end{align}
Combining \eqref{2-4}, \eqref{2-5}, \eqref{2-8} and \eqref{2-9}
yields
\begin{flalign}
&\sum_{n=0}^{\infty}(NB_2(1,5,5n+3)-NB_2(4,5,5n+3))q^n
=2\frac{(q^5;q^5)_{\infty}^3 }{ (q;q)_{\infty}} ,
 \label{3-12}\\[6pt]
&\sum_{n=0}^{\infty}(NB_2(2,5,5n+3)-NB_2(3,5,5n+3))q^n
=-\frac{(q^5;q^5)_{\infty}^3 }{ (q;q)_{\infty}}. \label{3-13}
\end{flalign}
Identity
  \eqref{1-8} follows from \eqref{3-12}
  and \eqref{3-13}. Furthermore,
   identities \eqref{1-9} and \eqref{1-10}
    follow from \eqref{2-15}, \eqref{3-12}
  and \eqref{3-13}.

Extracting those terms in which
 the power of $q$ is congruent to $4$
modulo $5$ in \eqref{3-3} and \eqref{3-4},
  then dividing by $q^4$
and replacing $q^5$ by $q$, we arrive at
\begin{align}
&\sum_{n=0}^{\infty}(NB_2(1,5,5n+4)-NB_2(4,5,5n+4))q^n
  \nonumber\\[6pt]
 =&V(q)X(4,1,0) +W(q) (X(4,2,1))
 \nonumber\\[6pt]
 =&-\frac{(q^5;q^5)_{\infty}^2}{
 (q^2,q^3;q^5)_{\infty}^2}
  \qquad ({\rm by} \  \eqref{2-3}\ {\rm
  and} \ \eqref{2-8}) \label{3-14}
\end{align}
and
\begin{align}
&\sum_{n=0}^{\infty}(NB_2(2,5,5n+4)-NB_2(3,5,5n+4))q^n
  \nonumber\\[6pt]
 =&V(q) X(2,2,0)
 +W(q) (X(3,1,0)-X(4,2,1))
   \nonumber \\[6pt]
 =&3\frac{(q^5;q^5)_{\infty}^2}
 {(q^2,q^3;q^5)_{\infty}^2},  \label{3-15}
  \qquad \qquad({\rm by} \  \eqref{2-4}\
  \eqref{2-7} \ {\rm
  and} \ \eqref{2-8})
\end{align}
from which with \eqref{3-14}, \eqref{3-15}
 identity \eqref{1-11}  follows.
  This completes the proof of Theorem
  \ref{Th1-1}.
 \qed

\section{Proof of Theorem \ref{Th1-2} }

Setting $k=3$ in \eqref{2-19} and utilizing
  \eqref{2-15}, \eqref{3-1}
    and \eqref{3-2}, we obtain
 \begin{align}
&\sum_{n\geq0}[NB_{3}(1,5,n)-NB_{3}(4,5,n)+2NB_{3}(2,5,n)-2NB_{3}(3,5,n)]q^n
\nonumber\\[6pt]
=&\frac{(q^{25};q^{25})_\infty^5}{(q^5;q^5)_\infty^6}
\bigg(\frac{V^5(q^5)}{W^4(q^5)}+3 q
\frac{V^4(q^5)}{W^3(q^5)}+4q^2\frac{V^3(q^5) }{W^2(q^5) }+7 q^3
\frac{V^2(q^5)}{W(q^5)}+11 q^4 V(q^5) +7 q^5 W(q^5) \nonumber \\[6pt]
&-4 q^6 \frac{W^2(q^5)}{V(q^5)} +3 q^7 \frac{W^3(q^5) }{V^2(q^5)
}-q^8 \frac{W^4(q^5)}{V^3(q^5)} +2q^9 \frac{W^5(q^5)
}{V^4(q^5)}\bigg)(R_1(q)-R_4(q))
\nonumber\\[6pt]
&+\frac{(q^{25};q^{25})_\infty^5
 }{(q^5;q^5)_\infty^6
  }\bigg(2 \frac{V^5(q^5)}{W^4(q^5)}
  +q \frac{V^4(q^5)}{W^3(q^5)}+3 q^2
\frac{V^3(q^5)}{W^2(q^5)}+4 q^3 \frac{V^2(q^5) }{W(q^5) }+7 q^4
V(q^5)-11 q^5 W(q^5) \nonumber\\[6pt]
 &+7 q^6 \frac{W^2(q^5)}{V(q^5)}
  -4 q^7
\frac{W^3(q^5)}{V^2(q^5)}+3 q^8 \frac{W^4(q^5) }{V^3(q^5) }-q^9
\frac{W^5(q^5)}{V^4(q^5)}\bigg)( R_2(q)-R_3(q)). \label{4-1}
\end{align}
Extracting
 those  terms in which the power of $q$ is
  congruent to 0 modulo 5 in
  \eqref{4-1}, then replacing $q^5$ by $q$,
we arrive at
\begin{flalign}
&\sum_{n\geq0}[NB_{3}(1,5,5n)-NB_{3}(4,5,5n)+2NB_{3}(2,5,5n)-2NB_{3}(3,5,5n)]q^n
  \nonumber\\[6pt]
=&\frac{(q^5;q^5)_\infty^5}{(q;q)_\infty^6}
 \bigg(\frac{V^5}{W^4}Y(1)
+3q\frac{V^4}{W^3}X(4,1,0) +4q\frac{V^3}{W^2}X(3,1,0)
+7q\frac{V^2}{W}X(2,1,0)+11qVX(1,1,0)
  \nonumber\\[6pt]
&+7qWY(1) -4q^2\frac{W^2}{V}X(4,1,0) +3q^2\frac{W^3}{V^2}X(3,1,0)
-q^2 \frac{W^4}{V^3}X(2,1,0) +2q^2 \frac{W^5}{V^4}X(1,1,0)
  \nonumber\\[6pt]
& +2 \frac{V^5}{W^4}Y(2) + q\frac{V^4}{W^3}X(2,2,0) +3q
\frac{V^3}{W^2}X(4,2,1) +4q \frac{V^2}{W}X(1,2,0) +7q V X(3,2,1)
 -11q W Y(2) \nonumber\\[6pt]
& +7q^2 \frac{W^2}{V}X(2,2,0) -4q^2 \frac{W^3}{V^2}X(4,2,1) +3q^2
\frac{W^4}{V^3}X(1,2,0)
 -q^2 \frac{W^5}{V^4}X(3,2,1) \bigg),
 \label{4-2}
\end{flalign}
where $V$ and $W$ are defined by \eqref{2-16}.
 Substituting
\eqref{2-3}--\eqref{2-12} into \eqref{4-2}
 yields
\begin{align}
&\sum_{n\geq0}[NB_{3}(1,5,5n)-NB_{3}(4,5,5n)+2NB_{3}(2,5,5n)-2NB_{3}(3,5,5n)]q^n
\nonumber\\[6pt]
=&22q^2\frac{(q,q^4;q^5)_{\infty}^2(q^5;q^5)_{\infty}^7}{(q;q)_{\infty}^6(q^2,q^3;q^5)_{\infty}^4}
+\frac{49q(q^2,q^3;q^5)_{\infty}(q^5;q^5)_{\infty}^7}{2(q;q)_{\infty}^6(q,q^4;q^5)_{\infty}^3}
+\frac{(q^2,q^3;q^5)_{\infty}^6
(q^5;q^5)_{\infty}^7}{2(q;q)_{\infty}^6(q,q^4;q^5)_{\infty}^8}
  \nonumber\\[6pt]
&-q\frac{(q^5;q^5)_{\infty}^5}{(q;q)_{\infty}^6(q^2,q^3;q^5)_{\infty}}
 -\frac{(q^2,q^3;q^5)_{\infty}^4(q^5;
 q^5)_{\infty}^5}{2(q;q)_{\infty}^6(
 q,q^4;q^5)_{\infty}^5}\nonumber\\[6pt]
  =&\frac{q(q^2,q^3;q^5)_{\infty}
  (q^5;q^5)_{\infty}^7}{2(q;q)_{\infty}^6
  (q,q^4;q^5)_{\infty}^3}\bigg(
  \frac{1}{R(q)^5}-11-R(q)^5+60
  +45R(q)^5\bigg)
  \nonumber \\[6pt]
&-\frac{(q^2,q^3;q^5)_{\infty} (q^5;q^5)_{\infty}^5}{2(q;q)_{\infty
}^6(q,q^4;q^5)_{\infty}^2}\bigg(
\frac{q^{\frac{3}{5}}}{R(q)^3}-3q^{\frac{3}{5}}R(q)^2+5q\frac{(q,q^4;q^5)_{\infty}^2}{(q^2,q^3;q^5)_{\infty}^2}\bigg)
 . \label{4-3}
\end{align}
Substituting
 \eqref{2-13} and
   \eqref{2-14} into
   \eqref{4-3}
    and using \eqref{2-18}, we arrive
at \eqref{1-12}.
 This completes the proof of Theorem \ref{Th1-2}.
  \qed

\section{Proof of Theorem \ref{Th1-3}}

Setting $k=4$ in \eqref{2-19} and employing
  \eqref{2-15}, \eqref{3-1}
    and \eqref{3-2}, we obtain
\begin{flalign}
&\sum_{n\geq0}
[NB_4(1,5,n)-NB_4(4,5,n)+2NB_4(2,5,n)-2NB_4(3,5,n)]q^n
  \nonumber \\[6pt]
=&\frac{(q^{25};q^{25})_\infty^{10}}{(q^5;q^5)_\infty^{12}}
\bigg(2q^{17}\frac{W^9(q^5)}{V^8(q^5)}
-3q^{16}\frac{W^8(q^5)}{V^7(q^5)} +8q^{15}\frac{W^7(q^5)}{V^6(q^5)}
-15q^{14}\frac{W^6(q^5)}{V^5(q^5)}
+30q^{13}\frac{W^5(q^5)}{V^4(q^5)}
\nonumber \\[6pt]
&-12q^{12}\frac{W^4(q^5)}{V^3(q^5)}
+38q^{11}\frac{W^3(q^5)}{V^2(q^5)} -13q^{10}\frac{W^2(q^5)}{V(q^5)}
+10q^9 W(q^5) +55q^8V(q^5) +74q^7\frac{V^2(q^5)}{W(q^5)}
\nonumber \\[6pt]
&+59q^6\frac{V^3(q^5)}{W^2(q^5)}+56q^5\frac{V^4(q^5)}{W^3(q^5)}
+40q^4\frac{V^5(q^5)}{W^4(q^5)}+20q^3 \frac{V^6(q^5)}{W^5(q^5)}
+9q^2\frac{V^7(q^5)}{W^6(q^5)}+4q\frac{V^8(q^5)}{W^7(q^5)}
\nonumber \\
&+\frac{V^9(q^5)}{W^8(q^5)}\bigg)(R_1(q)-R_4(q))
+\frac{(q^{25};q^{25})_\infty^{10}}{(q^5;q^5)_\infty^{12}}
\bigg(-q^{17}\frac{W^9(q^5)}{V^8(q^5)}+4q^{16}\frac{W^8(q^5)}{V^7(q^5)}-9
q^{15}\frac{W^7(q^5)}{V^6(q^5)}
\nonumber \\[6pt]
&+20q^{14}\frac{W^6(q^5)}{V^5(q^5)}
-40q^{13}\frac{W^5(q^5)}{V^4(q^5)}+56 q^{12}\frac{W^4(q^5)
}{V^3(q^5)}-59 q^{11}\frac{W^3(q^5)}{V^2(q^5)}+74q^{10}
\frac{W^2(q^5)}{V(q^5)}\nonumber \\[6pt]
 &-55 q^{9} W(q^5)+10 q^8 V(q^5)+13
q^7\frac{V^2(q^5)}{W(q^5)} + 38 q^6 \frac{V^3(q^5)}{W^2(q^5)}+12 q^5
\frac{V^4(q^5)}{W^3(q^5)}+12 q^5 \frac{V^4(q^5)}{W^3(q^5)}
\nonumber \\[6pt]
&+30 q^4 \frac{V^5(q^5)}{W^4(q^5)} +15 q^3
\frac{V^6(q^5)}{W^5(q^5)}+8 q^2 \frac{V^7(q^5) }{W^6(q^5)}+3 q
\frac{V^8(q^5)}{W^7(q^5)}
 +2\frac{V^9(q^5)}{W^8(q^5)}\bigg)(
 R_2(q)-R_3(q)). \label{5-1}
\end{flalign}
If we pick out
 those  terms in which the power of
  $q$ is
  congruent to 0 modulo 5 in \eqref{5-1}, then
   replace $q^5$ by $q$, we obtain
\begin{flalign}
&\sum_{n\geq0}[NB_{4}(1,5,5n)-NB_{4}(4,5,5n)+2NB_{4}(2,5,5n)-2NB_{4}(3,5,5n)]q^n
 \nonumber\\[6pt]
=&\frac{(q^5;q^5)_\infty^{10}}{(q;q)_\infty^{12}
 }\bigg(2q^4\frac{W^9}{V^8}X(3,1,0)
 -3q^4\frac{W^8}{V^7}X(4,1,0)
+8q^3\frac{W^7}{V^6}Y(1) -15q^3\frac{W^6}{V^5}X(1,1,0)
\nonumber\\[6pt]
&+30q^3\frac{W^5}{V^4}X(2,1,0) -12q^3\frac{W^4}{V^3}X(3,1,0)
+38q^3\frac{W^3}{V^2}X(4,1,0) -13q^2\frac{W^2}{V}Y(1)
\nonumber\\[6pt]
&+10q^2WX(1,1,0) +55q^2VX(2,1,0) +74q^2\frac{V^2}{W}X(3,1,0)
+59q^2\frac{V^3}{W^2}X(4,1,0)
\nonumber\\[6pt]
&+56q\frac{V^4}{W^3}Y(1) +40q\frac{V^5}{W^4}X(1,1,0)
+20q\frac{V^6}{W^5}X(2,1,0) +9q\frac{V^7}{W^6}X(3,1,0)
\nonumber\\[6pt]
&+4q\frac{V^8}{W^7}X(4,1,0) +\frac{V^9}{W^8}Y(1)
-q^4\frac{W^9}{V^8}X(4,2,1) +4q^4\frac{W^8}{V^7}X(2,2,0)
\nonumber\\[6pt]
&-9q^3\frac{W^7}{V^6}Y(2) +20q^3\frac{W^6}{V^5}X(3,2,1)
-40q^3\frac{W^5}{V^4}X(1,2,0) +56q^3\frac{W^4}{V^3}X(4,2,1)
\nonumber\\[6pt]
&-59q^3\frac{W^3}{V^2}X(2,2,0) +74q^2\frac{W^2}{V}Y(2)
-55q^2WX(3,2,1) +10q^2VX(1,2,0)
\nonumber\\[6pt]
&+13q^2\frac{V^2}{W}X(4,2,1) +38q^2\frac{V^3}{W^2}X(2,2,0)
+12q\frac{V^4}{W^3}Y(2) +30q\frac{V^5}{W^4}X(3,2,1)
\nonumber\\[6pt]
&+15q\frac{V^6}{W^5}X(1,2,0) +8q\frac{V^7}{W^6}X(4,2,1)
+3q\frac{V^8}{W^7}X(2,2,0)
 +2\frac{V^9}{W^8}Y(2)\bigg),
\label{5-2}
\end{flalign}
where $V$ and $W$ are defined by \eqref{2-16}.
 Substituting \eqref{2-3}--\eqref{2-12} into  \eqref{5-2} yields
\begin{flalign}
&\sum_{n\geq0}[NB_4(1,5,5n)-NB_4(4,5,5n)+2N
B_4(2,5,5n)-2NB_4(3,5,5n)]q^n
\nonumber\\[6pt]
 =&\frac{21}{2}q^4\frac{(q,q^4;q^5)_{\infty}^8(q^5;q^5)_{\infty}^{12}}{(q^2,q^3;q^5)_{\infty}^{10}(q;q)_{\infty}^{12}}
-174q^3\frac{(q,q^4;q^5)_{\infty}^3(q^5;q^5)_{\infty}^{12}}{(q^2,q^3;q^5)_{\infty}^5(q;q)_{\infty}^{12}}
+\frac{359}{2}q^2\frac{(q^5;q^5)_{\infty}^{12}}{(q,q^4;q^5)_{\infty}^2(q;q)_{\infty}^{12}}
\nonumber\\[6pt]
&+\frac{193}{2}q\frac{(q^2,q^3;q^5)_{\infty}^5(q^5;q^5)_{\infty}^{12}}{(q,q^4;q^5)_{\infty}^7(q;q)_{\infty}^{12}}
+\frac{1}{2}\frac{(q^2,q^3;q^5)_{\infty}^{10}(q^5;q^5)_{\infty}^{12}}{(q,q^4;q^5)_{\infty}^{12}(q;q)_{\infty}^{12}}
-\frac{3}{2}q^3\frac{(q,q^4;q^5)_{\infty}^6(q^5;q^5)_{\infty}^{10}}{(q^2,q^3;q^5)_{\infty}^7(q;q)_{\infty}^{12}}
\nonumber\\[6pt]
&-\frac{7}{2}q^2\frac{(q,q^4;q^5)_{\infty}(q^5;q^5)_{\infty}^{10}}{(q^2,q^3;q^5)_{\infty}^2(q;q)_{\infty}^{12}}
-18q\frac{(q^2,q^3;q^5)_{\infty}^3(q^5;q^5)_{\infty}^{10}}{(q,q^4;q^5)_{\infty}^4(q;q)_{\infty}^{12}}
-\frac{1}{2}\frac{(q^2,q^3;q^5 )_{\infty}^8(q^5;q^5)_{\infty}^{10}
}{(q,q^4;q^5)_{\infty}^9(q;q)_{\infty}^{12}} \nonumber
\\[6pt]
=&    \left(\frac{q^{\frac{3}{5}}}{R^3(q)}
-3q^{\frac{3}{5}}R^2(q)\right)\bigg( -\frac{7}{2}q^3\frac{(q,q^4;q^5
)_{\infty}^{6}(q^5;q^5)_{\infty}^{12}}
{(q^2,q^3;q^5)_{\infty}^{8}(q;q)_{\infty}^{12}}
 +6q^2\frac{(q,q^4;q^5)_{\infty} (q^5;q^5)_{\infty}^{12}}{(q^2,q^3;q^5
)_{\infty}^{3}(q;q)_{\infty}^{12}} \nonumber\\[6pt]
&-2q\frac{(q^2,q^3;q^5)_{\infty}^2
  (q^5;q^5)_{\infty}^{12}}{
  (q,q^4;q^5)_{\infty}^4(q;q)_{\infty}^{12}}
  -\frac{9}{2}\frac{(q^2,q^3;q^5)_{\infty}^7
(q^5;q^5)_{\infty}^{12}}{(q,q^4;q^5)_{\infty}^9(q;q)_{\infty}^{12}}
  +8q^2\frac{(q,q^4;q^5)_{\infty}^{4}(q^5;q^5)_{\infty}^{10}
}{(q^2,q^3;q^5)_{\infty}^{5}(q;q)_{\infty}^{12}} \nonumber\\[6pt]
&+\frac{q(q^5;q^5)_{\infty}^{10}}{2(q,q^4;q^5)_{\infty}(q;q)_{\infty}^{12}}
 +\frac{9}{2}\frac{(q^2,q^3;q^5)_{\infty}^5(q^5;q^5)_{\infty}^{10}}
{(q,q^4;q^5)_{\infty}^6(q;q)_{\infty}^{12}}\bigg)
-\frac{305}{2}q^3\frac{
(q,q^4;q^5)_{\infty}^3(q^5;q^5)_{\infty}^{12}}
{(q^2,q^3;q^5)_{\infty}^5(q;q)_{\infty}^{12}} \nonumber\\[6pt]
&+\frac{45}{2}q^3
\frac{(q,q^4;q^5)_{\infty}^6(q^5;q^5)_{\infty}^{10} }{
(q^2,q^3;q^5)_{\infty}^7(q;q)_{\infty}^{12}} +\frac{335}{2}q^2
\frac{ \left(q^5 ; q^5\right)_{\infty}^{12}}{\left(q, q^4 ;
 q^5\right)_{\infty}^2(q;q)_{\infty}^{12}}
-10q^2\frac{(q,q^4;q^5)_{\infty}(q^5;q^5)_{\infty}^{10}
 }{(q^2,q^3;q^5)_{\infty}^2(q;q)_{\infty}^{12}}
\nonumber\\[6pt]
 &+85q \frac{(q^2, q^3 ; q^5)_{\infty}^5\left(q^5;
q^5\right)_{\infty}^{12}}{\left(q, q^4 ;
q^5\right)_{\infty}^7(q;q)_{\infty}^{12}}
-5q\frac{(q^2,q^3;q^5)_{\infty}^3 (q^5;q^5)_{\infty}^{10}}{
(q,q^4;q^5)_{\infty}^{4}(q;q)_{\infty}^{12} }
+5\frac{(q^2,q^3;q^5)_{\infty}^{10}
 (q^5;q^5)_{\infty}^{12}}{
(q,q^4;q^5)_{ \infty}^{12} (q;q)_{\infty}^{12}}
\nonumber\\[6pt]
 &-5\frac{(q^2,q^3;q^5)_{\infty}^8
 (q^5;q^5)_{\infty}^{10}}{(q,q^4;
 q^5)_{\infty}^{9}(q;q)_{\infty}^{12}},
  \label{5-3}
\end{flalign}
 from which with
  \eqref{2-14},
 identity \eqref{1-13} follows.

Extracting those terms in which the power of $q$ is congruent to $3$
modulo $5$ in \eqref{5-1}, then dividing by $q^3$ and replacing
$q^5$ by $q$, we obtain
\begin{flalign}
&\sum_{n\geq0}[NB_{4}(1,5,5n+3)-NB_{4}(4,5,5n+3)+2NB_{4}(2,5,5n+3)-2NB_{4}(3,5,5n+3)]q^n
\nonumber\\[6pt]
 =&\frac{(q^5;q^5)_\infty^{10}}{(q;q)_\infty^{12}
 }\bigg(2q^3\frac{W^9}{V^8}X(1,1,0)
 -3q^3\frac{W^8}{V^7}X(2,1,0)
+8q^3\frac{W^7}{V^6}X(3,1,0) -15q^3\frac{W^6}{V^5}X(4,1,0)
\nonumber\\[6pt]
&+30q^2\frac{W^5}{V^4}Y(1) -12q^2\frac{W^4}{V^3}X(1,1,0)
+38q^2\frac{W^3}{V^2}X(2,1,0) -13q^2\frac{W^2}{V}X(3,1,0)
\nonumber\\[6pt]
&+10q^2WX(4,1,0) +55qVY(1) +74q\frac{V^2}{W}X(1,1,0)
+59q\frac{V^3}{W^2}X(2,1,0)
\nonumber\\[6pt]
&+56q\frac{V^4}{W^3}X(3,1,0) +40q\frac{V^5}{W^4}X(4,1,0)
+20\frac{V^6}{W^5}Y(1) +9\frac{V^7}{W^6}X(1,1,0)
\nonumber\\[6pt]
&+4\frac{V^8}{W^7}X(2,1,0) +\frac{V^9}{W^8}X(3,1,0)
-q^3\frac{W^9}{V^8}X(3,2,1) +4q^3\frac{W^8}{V^7}X(1,2,0)
\nonumber\\[6pt]
&-9q^3\frac{W^7}{V^6}X(4,2,1) +20q^3\frac{W^6}{V^5}X(2,2,0)
-40q^2\frac{W^5}{V^4}Y(2) +56q^2\frac{W^4}{V^3}X(3,2,1)
\nonumber\\[6pt]
&-59q^2\frac{W^3}{V^2}X(1,2,0) +74q^2\frac{W^2}{V}X(4,2,1)
-55q^2WX(2,2,0) +10qVY(2)
\nonumber\\[6pt]
&+13q\frac{V^2}{W}X(3,2,1) +38q\frac{V^3}{W^2}X(1,2,0)
+12q\frac{V^4}{W^3}X(4,2,1) +30q\frac{V^5}{W^4}X(2,2,0)
\nonumber\\[6pt]
&+15\frac{V^6}{W^5}Y(2) +8\frac{V^7}{W^6}X(3,2,1)
+3\frac{V^8}{W^7}X(1,2,0) +2\frac{V^9}{W^8} X(4,2,1)\bigg).
\label{5-4}
\end{flalign}
Substituting \eqref{2-3}--\eqref{2-12} into \eqref{5-4}, we  arrive
at   \eqref{1-14}.

Extracting those terms in which the power of $q$ is congruent to $4$
modulo $5$ in \eqref{5-1}, then dividing by $q^4$ and replacing
$q^5$ by $q$, we obtain
\begin{flalign}\label{5-5}
&\sum_{n\geq0}[NB_{4}(1,5,5n+4)-NB_{4}(4,5,5n+4)+2NB_{4}(2,5,5n+4)-2NB_{4}(3,5,5n+4)]q^n
\nonumber\\[6pt]
 =&\frac{(q^5;q^5)_\infty^{10}}{(q;q)_\infty^{12}
 }\bigg(2q^3\frac{W^9}{V^8}X(2,1,0)
 -3q^3\frac{W^8}{V^7}X(3,1,0)
+8q^3\frac{W^7}{V^6}X(4,1,0) -15q^2\frac{W^6}{V^5}Y(1)
\nonumber\\[6pt]
&+30q^2\frac{W^5}{V^4}X(1,1,0) -12q^2\frac{W^4}{V^3}X(2,1,0)
+38q^2\frac{W^3}{V^2}X(3,1,0) -13q^2\frac{W^2}{V}X(4,1,0)
\nonumber\\[6pt]
&+10qWY(1)+55qVX(1,1,0) +74q\frac{V^2}{W}X(2,1,0)
+59q\frac{V^3}{W^2}X(3,1,0)
\nonumber\\[6pt]
&+56q\frac{V^4}{W^3}X(4,1,0) +40\frac{V^5}{W^4}Y(1)
+20\frac{V^6}{W^5}X(1,1,0) +9\frac{V^7}{W^6}X(2,1,0)
\nonumber\\[6pt]
&+4\frac{V^8}{W^7}X(3,1,0) +\frac{V^9}{W^8}X(4,1,0)
-q^3\frac{W^9}{V^8}X(1,2,0) +4q^3\frac{W^8}{V^7}X(4,2,1)
\nonumber\\[6pt]
&-9q^3\frac{W^7}{V^6}X(2,2,0) +20q^2\frac{W^6}{V^5}Y(2)
-40q^2\frac{W^5}{V^4}X(3,2,1) +56q^2\frac{W^4}{V^3}X(1,2,0)
\nonumber\\[6pt]
&-59q^2\frac{W^3}{V^2}X(4,2,1) +74q^2\frac{W^2}{V}X(2,2,0) -55qWY(2)
+10qVX(3,2,1)
\nonumber\\[6pt]
&+13q\frac{V^2}{W}X(1,2,0) +38q\frac{V^3}{W^2}X(4,2,1)
+12q\frac{V^4}{W^3}X(2,2,0) +30\frac{V^5}{W^4}Y(2)
\nonumber\\[6pt]
&+15\frac{V^6}{W^5}X(3,2,1) +8\frac{V^7}{W^6}X(1,2,0)
+3\frac{V^8}{W^7}X(4,2,1) +2\frac{V^9}{W^8}X(2,2,0)\bigg).
\end{flalign}
Substituting \eqref{2-3} -- \eqref{2-13}
 into \eqref{5-5}, we can
get   \eqref{1-15}.  This completes the proof of Theorem
\ref{Th1-3}.
 \qed

\section{Proof of Theorem \ref{Th1-4}}

Setting $k=5$ in \eqref{2-19} and
 utilizing
  \eqref{2-15}, \eqref{3-1}
    and \eqref{3-2}, we obtain
\begin{flalign}
&\sum_{n\geq 0} [N B_5(1,5, n)-N B_5(4,5, n)+3 NB_5(2,5, n)-3 N
B_5(3,5, n)]q^n \nonumber\\[6pt]
 =&\frac{(q^{25};q^{25})_\infty^{15}}{(q^5;q^5
  )_\infty^{18}
 }\bigg(3q^{25}\frac{W^{13}(q^5)}{V^{12}(q^5)
 }-8q^{24}\frac{W^{12}(q^5)}{V^{11}(q^5)}+
 24q^{23}\frac{W^{11}(q^5)}{V^{10}(q^5)}
-57q^{22}\frac{W^{10}(q^5)}{V^9(q^5)}
 +131q^{21}
\frac{W^9(q^5)}{V^8(q^5)} \nonumber\\[6pt]
&-183 q^{20}\frac{W^8(q^5)}{V^7(q^5)}+324q^{19}
\frac{W^7(q^5)}{V^6(q^5)} -397q^{18} \frac{W^6(q^5)}{V^5(q^5)}
+471q^{17} \frac{W^5(q^5)}{V^4(q^5)} -243q^{16}
 \frac{W^4(q^5)}{V^3(q^5)}
  \nonumber\\[6pt]
& +504q^{15}
 \frac{W^3(q^5)}{V^2(q^5)}+48q^{14} \frac{W^2(q^5)}{V(q^5)}+156 q^{13} W(q^5)+242q^{12}
V(q^5) +714 q^{11} \frac{V^2(q^5)
}{W(q^5)}+678 q^{10} \frac{V^3(q^5)}{W^2(q^5)} \nonumber\\[6pt]
&+801q^9 \frac{V^4(q^5)}{W^3(q^5)}
 +747q^8 \frac{V^5(q^5)}{W^4(q^5)}+579 q^7
\frac{V^6(q^5)}{W^5(q^5)}+368 q^6
 \frac{V^7(q^5)}{W^6(q^5)}+231 q^5
\frac{V^8(q^5)}{W^7(q^5)} \nonumber\\[6pt]
&+117q^4\frac{V^9(q^5)}{W^8(q^5)} +49q^3\frac{V^{10}(q^5)}{W^9(q^5)}
+18q^2\frac{V^{11}(q^5)}{W^{10}(q^5)}
+6q\frac{V^{12}(q^5)}{W^{11}(q^5)} +\frac{V^{13}(q^5)}{W^{12}(q^5) }
\bigg)(R_1(q)-R_4(q))  \nonumber \\[6pt]
&+\frac{(q^{25};q^{25})_\infty^{15}}{(q^5;q^5
  )_\infty^{18}
 }\bigg(-2q^{25}\frac{W^{13}(q^5)}{V^{12}(q^5)}
 +9q^{24}\frac{W^{12}(q^5)}{V^{11}(q^5)}
 -27q^{23}\frac{W^{11}(q^5)}{V^{10}(q^5)}
+71q^{22}\frac{W^{10}(q^5)}{V^9(q^5)} -168q^{21}
\frac{W^9(q^5)}{V^8(q^5)} \nonumber\\[6pt]
 &+309 q^{20}
 \frac{W^8(q^5)}{V^7(q^5)}-502 q^{19}
\frac{W^7(q^5)}{V^6(q^5)} +756q^{18} \frac{W^6(q^5)}{V^5(q^5) }-963
q^{17} \frac{W^5(q^5)}{V^4(q^5)}+954 q^{16}
\frac{W^4(q^5)}{V^3(q^5)}
\nonumber\\[6pt]
&-897 q^{15} \frac{W^3(q^5)}{V^2(q^5)}
  +771 q^{14} \frac{W^2(q^5)}{V(q^5)} -313 q^{13} W(q^5)+99 q^{12}
V(q^5) -267 q^{11}\frac{V^2(q^5)}{W(q^5)} +351
 q^{10}
\frac{V^3(q^5)}{W^2(q^5)}\nonumber\\[6pt]
&+27 q^9 \frac{V^4(q^5)}{W^3(q^5)}
 +294q^8\frac{V^5(q^5)}{W^4(q^5)}+263 q^7
\frac{V^6(q^5)}{W^5(q^5)}+246 q^6 \frac{V^7(q^5)}{W^6(q^5)} +132 q^5
\frac{V^8(q^5)}{W^7(q^5)}  +109q^4\frac{V^9(q^5)}{W^8(q^5)
}\nonumber\\[6pt]
&+48q^3\frac{V^{10}(q^5)}{W^9(q^5)}+
21q^2\frac{V^{11}(q^5)}{W^{10}(q^5)}
+7q\frac{V^{12}(q^5)}{W^{11}(q^5)} +3\frac{V^{13}(q^5)}{W^{12}(q^5)}
\bigg)(R_2(q)-R_3(q)). \label{6-1}
\end{flalign}
Extracting those terms in which the power of $q$ is congruent to $2$
modulo $5$ in \eqref{6-1},
 then dividing by $q^2$ and replacing $q^5$
by $q$, we obtain
\begin{flalign}
&\sum_{n\geq 0} [N B_5(1,5, 5n+2)-N B_5(4,5, 5n+2)+3 NB_5(2,5,
5n+2)-3 N B_5(3,5, 5n+2)]q^n
\nonumber\\[6pt]
 =&\frac{(q^{5};q^{5})_\infty^{15}}{(q;q
  )_\infty^{18}
 }\bigg(3q^{5}\frac{W^{13}}{V^{12}}X(2,1,0)
-8q^{5}\frac{W^{12}}{V^{11}}X(3,1,0)
+24q^{5}\frac{W^{11}}{V^{10}}X(4,1,0) -57q^{4}\frac{W^{10}}{V^9}Y(1)
\nonumber\\[6pt]
&+131q^{4} \frac{W^9}{V^8}X(1,1,0) -183q^{4}\frac{W^8}{V^7}X(2,1,0)
+324q^{4}\frac{W^7}{V^6}X(3,1,0) -397q^{4}\frac{W^6}{V^5}X(4,1,0)
\nonumber\\[6pt]
&+471q^{3}\frac{W^5}{V^4}Y(1) -243q^{3}\frac{W^4}{V^3}X(1,1,0)
+504q^{3}\frac{W^3}{V^2}X(2,1,0) +48q^{3}\frac{W^2}{V}X(3,1,0)
\nonumber\\[6pt]
&+156q^{3}W X(4,1,0) +242q^{2}V Y(1) +714q^2\frac{V^2}{W}X(1,1,0)
+678q^2\frac{V^3}{W^2}X(2,1,0)
\nonumber\\[6pt]
&+801q^2\frac{V^4}{W^3}X(3,1,0) +747q^2\frac{V^5}{W^4}X(4,1,0)
+579q\frac{V^6}{W^5}Y(1)+368q\frac{V^7}{W^6}X(1,1,0)
\nonumber\\[6pt]
&+231q\frac{V^8}{W^7}X(2,1,0) +117q\frac{V^9}{W^8}X(3,1,0)
+49q\frac{V^{10}}{W^9}X(4,1,0) +18\frac{V^{11}}{W^{10}}Y(1)
\nonumber\\[6pt]
&+6\frac{V^{12}}{W^{11}}X(1,1,0) +\frac{V^{13}}{W^{12}}X(2,1,0)
-2q^{5}\frac{W^{13}}{V^{12}}X(1,2,0)
+9q^{5}\frac{W^{12}}{V^{11}}X(4,2,1)
\nonumber\\[6pt]
&-27q^{5}\frac{W^{11}}{V^{10}}X(2,2,0)
+71q^{4}\frac{W^{10}}{V^9}Y(2) -168q^{4}\frac{W^9}{V^8}X(3,2,1)
+309q^4\frac{W^8}{V^7}X(1,2,0)
\nonumber\\[6pt]
&-502 q^{4} \frac{W^7}{V^6}X(4,2,1) +756q^{4}\frac{W^6}{V^5}X(2,2,0)
-963 q^{3}\frac{W^5}{V^4}Y(2) +954q^{3}\frac{W^4}{V^3}X(3,2,1)
\nonumber\\[6pt]
&-897q^{3}\frac{W^3}{V^2}X(1,2,0) +771q^{3}\frac{W^2}{V}X(4,2,1)
-313q^{3}W X(2,2,0) +99q^{2}V Y(2)
\nonumber\\[6pt]
&-267 q^2\frac{V^2}{W}X(3,2,1) +351q^2\frac{V^3}{W^2}X(1,2,0) +27q^2
\frac{V^4}{W^3}X(4,2,1) +294q^2\frac{V^5}{W^4}X(2,2,0)
\nonumber\\[6pt]
&+263 q \frac{V^6}{W^5}Y(2) +246q\frac{V^7}{W^6}X(3,2,1) +132q
\frac{V^8}{W^7}X(1,2,0) +109q\frac{V^9}{W^8}X(4,2,1)
\nonumber\\[6pt]
&+48q\frac{V^{10}}{W^9}X(2,2,0) +21\frac{V^{11}}{W^{10}}Y(2)
+7\frac{V^{12}}{W^{11}}X(3,2,1) +3\frac{V^{13}}{W^{12}}X(1,2,0)
 \bigg). \label{6-2}
\end{flalign}
Substituting \eqref{2-3}--\eqref{2-12} into \eqref{6-2},
 we arrive at \eqref{1-16}.

  Extracting those terms in which the
power of $q$ is congruent to $4$ modulo $5$ in \eqref{6-1}, then
dividing by $q^4$ and replacing $q^5$ by $q$, we obtain
\begin{flalign}
&\sum_{n\geq0}[NB_5(1,5,5n+4)-NB_5(4,5,5n+4)+3NB_5(2,5,5n+4)-3NB_5(3,5,5n+4)]q^n
\nonumber\\[6pt]
 =&\frac{(q^5;q^5)_\infty^{15}
 }{(q;q)_\infty^{18}}
  \bigg(3q^{5}\frac{W^{13}}{V^{12}}X(4,1,0)
-8q^{4}\frac{W^{12}}{V^{11}}Y(1)
+24q^{4}\frac{W^{11}}{V^{10}}X(1,1,0)
-57q^{4}\frac{W^{10}}{V^9}X(2,1,0)
\nonumber\\[6pt]
&+131q^{4} \frac{W^9}{V^8}X(3,1,0) -183q^{4}\frac{W^8}{V^7}X(4,1,0)
+324q^{3}\frac{W^7}{V^6}Y(1) -397q^{3}\frac{W^6}{V^5}X(1,1,0)
\nonumber\\[6pt]
&+471q^{3} \frac{W^5}{V^4}X(2,1,0) -243q^{3}\frac{W^4}{V^3}X(3,1,0)
+504q^{3}\frac{W^3}{V^2}X(4,1,0) +48q^{2}\frac{W^2}{V}Y(1)
\nonumber\\[6pt]
&+156q^2 W X(1,1,0) +242q^2V X(2,1,0) +714q^2 \frac{V^2}{W}X(3,1,0)
+678q^2\frac{V^3}{W^2}X(4,1,0)
\nonumber\\[6pt]
&+801q\frac{V^4}{W^3}Y(1) +747q\frac{V^5}{W^4}X(1,1,0)
+579q\frac{V^6}{W^5}X(2,1,0) +368q\frac{V^7}{W^6}X(3,1,0)
\nonumber\\[6pt]
&+231q\frac{V^8}{W^7}X(4,1,0) +117\frac{V^9}{W^8}Y(1)
+49\frac{V^{10}}{W^9}X(1,1,0) +18\frac{V^{11}}{W^{10}}X(2,1,0)
\nonumber\\[6pt]
&+6\frac{V^{12}}{W^{11}}X(3,1,0) +\frac{V^{13}}{W^{12}}X(4,1,0)
-2q^{5}\frac{W^{13}}{V^{12}}X(2,2,0)
+9q^{4}\frac{W^{12}}{V^{11}}Y(2)
\nonumber\\[6pt]
&-27q^{4}\frac{W^{11}}{V^{10}}X(3,2,1)
+71q^{4}\frac{W^{10}}{V^9}X(1,2,0) -168q^{4}\frac{W^9}{V^8}X(4,2,1)
+309q^{4}\frac{W^8}{V^7}X(2,2,0)
\nonumber\\[6pt]
&-502q^{3}\frac{W^7}{V^6}Y(2) +756q^{3}\frac{W^6}{V^5}X(3,2,1)
-963q^{3}\frac{W^5}{V^4}X(1,2,0) +954q^{3}\frac{W^4}{V^3}X(4,2,1)
\nonumber\\[6pt]
&-897 q^{3} \frac{W^3}{V^2}X(2,2,0) +771q^{2}\frac{W^2}{V}Y(2)
-313q^2WX(3,2,1) +99q^2VX(1,2,0)
\nonumber\\[6pt]
&-267q^2\frac{V^2}{W}X(4,2,1) +351q^2\frac{V^3}{W^2}X(2,2,0)
+27q\frac{V^4}{W^3}Y(2) +294q\frac{V^5}{W^4}X(3,2,1)
\nonumber\\[6pt]
&+263q\frac{V^6}{W^5}X(1,2,0) +246q\frac{V^7}{W^6}X(4,2,1)
+132q\frac{V^8}{W^7}X(2,2,0) +109\frac{V^9}{W^8}Y(2)
\nonumber\\[6pt]
&+48\frac{V^{10}}{W^9}X(3,2,1) +21\frac{V^{11}}{W^{10}}X(1,2,0)
+7\frac{V^{12}}{W^{11}}X(4,2,1) +3\frac{V^{13}}{W^{12}}
X(2,2,0)\bigg). \label{6-3}
\end{flalign}
Substituting \eqref{2-3}--\eqref{2-12} into \eqref{6-3} yields
\begin{flalign}
&\sum_{n\geq 0}[N B_5(1,5,5 n+4)-N B_5(4,5,5 n+4)+3 N B_5(2,5,5
n+4)-3N B_5(3,5,5 n+4)] q^n
  \nonumber \\[6pt]
 =&-\frac{11}{2}q^5\frac{(q,q^4;q^5)_{\infty}^{13}(q^5;q^5)_{\infty}^{17}}{(q^2,q^3;q^5)_{\infty}^{15}(q;q)_{\infty}^{18}}
+\frac{1655}{2}q^4\frac{(q,q^4;q^5)_{\infty}^{8}(q^5;q^5)_{\infty}^{17}}{(q^2,q^3;q^5)_{\infty}^{10}(q;q)_{\infty}^{18}}
-\frac{6325}{2}q^3 \frac{(q, q^4 ; q^5)_{\infty}^3(q^5
;q^5)_{\infty}^{17}}{(q^2, q^3 ; q^5)_{\infty}^5(q;q)_{\infty}^{18}}
\nonumber \\[6pt]
&+\frac{3985}{2}q^2 \frac{(q^5 ;q^5)_{\infty}^{17}}{(q, q^4;
q^5)_{\infty}^2(q;q)_{\infty}^{18}}
 +2065q
\frac{(q^2, q^3 ; q^5)_{\infty}^5(q^5 ; q^5)_{\infty}^{17}}{(q, q^4
; q^5)_{\infty}^7(q;q)_{\infty}^{18}} +136\frac{(q^2, q^3 ;
q^5)_{\infty}^{10}(q^5;q^5)_{\infty}^{17}}{(q, q^4 ;
q^5)_{\infty}^{12}(q;q)_{\infty}^{18}}
\nonumber \\[6pt]
&+\frac{3}{2}q^4\frac{(q,q^4;q^5)_{\infty}^{11}(q^5;q^5)_{\infty}^{15}}{(q^2,q^3;q^5)_{\infty}^{12}(q;q)_{\infty}^{18}}
-47q^3\frac{(q,q^4;q^5)_{\infty}^6(q^5;q^5)_{\infty}^{15}}{(q^2,q^3;q^5)_{\infty}^7(q;q)_{\infty}^{18}}
-\frac{183}{2}q^2
\frac{(q,q^4;q^5)_{\infty}(q^5;q^5)_{\infty}^{15}}{(q^2, q^3 ;
q^5)_{\infty}^2(q;q)_{\infty}^{18}}
\nonumber \\[6pt]
&-243q \frac{(q^2,q^3;q^5)_{\infty}^3(q^5;q^5)_{\infty}^{15}}{(q,
q^4 ; q^5)_{\infty}^4(q;q)_{\infty}^{18}} -46\frac{(q^2, q^3 ;
q^5)_{\infty}^{8}(q^5;q^5)_{\infty}^{15}}{(q, q^4 ;
q^5)_{\infty}^{9}(q;q)_{\infty}^{18}} \nonumber\\[6pt]
=&\left( \frac{q^{\frac{3}{5}}}{R^3(q)} -3q^{\frac{3}{5}}R^2(q)
\right)\bigg(q^4\frac{(q,q^4;q^5)_{\infty}^{11}
(q^5;q^5)_{\infty}^{17}}{ (q^2,q^3;q^5
)_{\infty}^{13}(q;q)_{\infty}^{18}} -\frac{q^3}{2} \frac{
(q,q^4;q^5)_{\infty}^{6} (q^5;q^5)_{\infty}^{17}}{
 (q^2,q^3;q^5)_{\infty}^{8}(q;q)_{\infty}^{18}}
\nonumber\\[6pt]
 &-\frac{q^3}{2} \frac{(q,q^4;q^5
)_{\infty}^{9}(q^5;q^5)_{\infty}^{15}}{(q^2,q^3;q^5)_{\infty}^{10}(q;q)_{\infty}^{18}}
+1054q^2\frac{(q,q^4;q^5)_{\infty}(q^5;q^5)_{\infty}^{17}}{(q^2,q^3;q^5)_{\infty}^3(q;q)_{\infty}^{18}}
+\frac{31q^2}{2}\frac{(q,q^4;q^5)_{\infty}^{4}(q^5;q^5)_{\infty}^{15}}{(q^2,q^3;q^5)_{\infty}^{5}(q;q)_{\infty}^{18}}
\nonumber\\[6pt]
&+\frac{3q(q^5;q^5)_{\infty}^{15}}{2(q,q^4;q^5
)_{\infty}(q;q)_{\infty}^{18}} -2q\frac{(q^2,q^3;q^5)_\infty^2
   (q^5;q^5)_\infty^{17} }{
    (q,q^4;q^5)_\infty^4(q;q)_{\infty}^{18}}
-\frac{3}{2}\frac{
 (q^2,q^3;q^5)_{\infty}^{7}(q^5;q^5)_{
 \infty}^{17}}{(q,q^4;q^5)_{\infty}^{9}(q;q)_{\infty}^{18}}
\nonumber\\[6pt]
&+\frac{3}{2}\frac{(q^2,q^3;
q^5)_{\infty}^{5}(q^5;q^5)_{\infty}^{15}}{(q,q^4;q^5
)_{\infty}^{6}(q;q)_{\infty}^{18}} \bigg)
-\frac{5}{2}q^5\frac{(q,q^4;q^5)_{\infty
}^{13}(q^5;q^5)_{\infty}^{18}}{
(q^2,q^3;q^5)_{\infty}^{15}(q;q)_{\infty}^{18}}
+825q^4\frac{(q,q^4;q^5)_{\infty}^{8} (q^5;q^5)_{\infty}^{17}}{
(q^2,q^3;q^5)_{\infty}^{10}(q;q)_{\infty}^{18}}
  \nonumber \\[6pt]
& -\frac{205q^2}{2}
 \frac{(q,q^4;q^5)_{\infty}(q^5;q^5)_{\infty}^{15}}
 {(q^2,q^3;q^5)_{\infty}^2(q;q)_{\infty}^{18}}
+\frac{1865q^2}{2}
\frac{(q^5;q^5)_{\infty}^{17}}{(q,q^4;q^5)_{\infty}^2(q;q)_{\infty}^{18}}
 -240q\frac{(q^2,q^3;q^5)_{\infty}^3(q^5;q^5)_{\infty}^{15}
}{(q,q^4;q^5)_{\infty}^4(q;q)_{\infty}^{18}}
  \nonumber \\[6pt]
&+\frac{4125q}{2}
\frac{(q^2,q^3;q^5)_{\infty}^5(q^5;q^5)_{\infty}^{17}}{(q,q^4;q^5)_{\infty}^7(q;q)_{\infty}^{18}}
 +\frac{275}{2}\frac{(q^2,q^3;q^5)_{\infty}^{10}(q^5;q^5)_{\infty}^{17}
}{(q,q^4;q^5)_{\infty}^{12}(q;q)_{\infty}^{18}}
-\frac{95}{2}\frac{(q^2,q^3;q^5)_{\infty}^{8}(q^5;q^5)_{\infty}^{15}}{(q,q^4;q^5)_{\infty}^{9}(q;q)_{\infty}^{18}}
 \label{6-4}
\end{flalign}
Identity \eqref{1-17} follows from
  \eqref{2-14} and \eqref{6-4}.
   This completes the proof.
    \qed

\section{Concluding remarks}

As seen in Introduction,  a number of nice
 results on Beck's
 partition   statistics
  $NT(r,m,n)$ and $M_{\omega}(r,m,n)$
   have been proved
 in recent years. Motivated by those
  work, we establish some identities on
   $NB_k(r,m,n)$ which is a
     partition  statistic
  of $k$-colored partitions. These identities
are analogous to
            Ramanujan's ``most beautiful
     identity". From these identities,
      one can easily deduce some congruences
       modulo 5 for  $NB_k(r,5,n)$
        with $2\leq k\leq 5$ proved by
        Lin, Peng and Toh \cite{Lin}.
 A natural question
 is to deduce   identities
 on  $NB_k(r,m,n)$  with $m\geq 7$
   which imply Lin, Peng and Toh's
      congruences
 modulo other moduli, such as 7, 11, 13.

\vspace{1cm}

 \noindent{\bf Acknowledgments.}
 This work was supported by
 the National Science Foundation
 of China (grant No.~11971203) and
     the Natural Science Foundation of
   Jiangsu Province of China (BK20221383).

   \noindent{\bf Declaration of
  Competing Interest.}
 The authors  declared that  they have
  no
 conflicts of interest to this work.

   \noindent{\bf Data Availability
   Statements.} Data
sharing not applicable to this article as no datasets were generated
or
 analysed during the current study.

\end{document}